\documentclass[11pt]{article}
\usepackage{amsmath,amssymb}

\newtheorem{theo}{Theorem}[section]

\newtheorem{lemma}[theo]{Lemma}
\newtheorem{claim}[theo]{Claim}
\newtheorem{coro}[theo]{Corollary}
\newcommand\npf{\mbox{ }\hfill\sqr\vskip6pt}
\def\sqr{$\vcenter{\hrule height.2mm
\hbox{\vrule width.2mm height2mm\kern2mm
\vrule width.2mm}\hrule height.2mm}$}

\begin{document}

\title{Equitable coloring of $k$-uniform hypergraphs}

\author{Raphael Yuster\\
Department of Mathematics\\ University of Haifa at Oranim\\ Tivon 36006,
Israel.\\ e-mail: raphy@research.haifa.ac.il}

\date{}

\maketitle

\begin{abstract}
Let $H$ be a $k$-uniform hypergraph with $n$ vertices. A {\em strong $r$-coloring} is a partition of the vertices
into $r$ parts, such that each edge of $H$ intersects each part. A strong $r$-coloring
is called {\em equitable} if the size of each part is $\lceil n/r \rceil$ or $\lfloor n/r \rfloor$.
We prove that for all $a \geq 1$, if the maximum degree of $H$ satisfies $\Delta(H) \leq k^a$
then $H$ has an equitable coloring with $\frac{k}{a \ln k}(1-o_k(1))$ parts.
In particular, every $k$-uniform hypergraph with maximum degree $O(k)$ has an equitable coloring
with $\frac{k}{\ln k}(1-o_k(1))$ parts.
The result is asymptotically tight. The proof uses a double application of the non-symmetric
version of the Lov\'asz Local Lemma.\\
{\bf Ams classification code:} 05C15
\end{abstract}

\setcounter{page}{1}

\section{Introduction}
All hypergraphs considered here are finite. For standard
terminology the reader is referred to \cite{BoMu}. Let $H$ be a $k$-uniform hypergraph with $n$ vertices.
A {\em strong $r$-coloring} is a partition of the vertices of $H$
into $r$ parts, such that each edge of $H$ intersects each part.
(A {\em weak} $r$-coloring is a coloring where no edge is monochromatic.)
A strong $r$-coloring is called {\em equitable} if the size of each part is $\lceil n/r \rceil$  or $\lfloor n/r \rfloor$.
Let $c(H)$ denote the maximum possible number of parts in a strong coloring of $H$.
Let $ec(H)$ denote the maximum possible number of parts in an equitable
coloring of $H$. Trivially, $1 \leq ec(H) \leq c(H) \leq k$. In general, $k$ could be large
and still $ec(H)=c(H)=1$, if we do not impose upper bounds on the maximum degree.
Consider the complete $k$-uniform hypergraph on $2k$ vertices.
Trivially, it has $c(H)=1$, and the maximum degree is less than $4^k$.
In this paper we prove that $c(H)$ and $ec(H)$ are quite large if the maximum degree
is bounded by a polynomial in $k$. In fact, we get the following tight result:
\begin{theo}
\label{t1}
Let $a \geq 1$. Let $H$ be a $k$-uniform hypergraph with maximum degree at most $k^a$.
Then, $ec(H) \geq \frac{k}{a \ln k}(1-o_k(1)).$
The result is tight. For all $a \geq 1$, there exist $k$-uniform hypergraphs $H$ with maximum degree at most $k^a$
and $c(H)\leq \frac{k}{a\ln k}(1+o_k(1)).$
\end{theo}
The tightness is shown by a construction of a random hypergraph with appropriate
parameters. Alon \cite{Al} has shown that there exist $k$-uniform hypergraphs with
$n$ vertices and maximum degree at most $k$ that do not have a vertex cover of size less than
$(n \ln k /k)(1-o_k(1))$. In particular, no strong coloring (moreover an equitable one) could have more than
$(k/\ln k)(1+o_k(1))$ parts. For completeness, we show a general construction valid for
all $a \geq 1$ in Section 3. The proof of the main result appears in Section 2. The final section
contains some concluding remarks.

\section{Proof of the main result}
In the proof of Theorem \ref{t1} we need to use the
Lov\'asz Local Lemma \cite{ErLo} in
its strongest form, known as the {\em nonsymmetric version}. Here it is,
following the notations in \cite{AlSp} (which also contains a simple proof of
the lemma). Let $A_1,\ldots,A_n$ be events in an arbitrary probability space.
A directed graph $D=(V,E)$ on the set of vertices $V=[n]$ is called a {\em
dependency digraph} for the events $A_1,\ldots,A_n$ if for each $i$,
$i=1,\ldots,n$, the event $A_i$ is mutually independent of all the events
$\{A_j ~ : ~ (i,j) \notin E\}$.
\begin{lemma}[The Local Lemma, nonsymmetric version]
\label{l21}
Let $A_1,\ldots,A_n$ be events in an arbitrary probability space and let
$D=(V,E)$ be a corresponding dependency digraph. If $x_1,\ldots,x_n$ are real
numbers such that $0 \leq x_i < 1$ and $\Pr[A_i] \leq x_i \prod_{(i,j) \in
E}(1-x_j)$ for all $i=1,\ldots,n$, then, with positive probability no event
$A_i$ holds. \npf
\end{lemma}
If the maximum outdegree in $D$ is at most $d \geq 1$ and each $A_i$ has $\Pr[A_i] \leq
p$ then, by assigning $x_i=1/(d+1)$ we immediately get:
\begin{coro}[The Local Lemma, symmetric version]
\label{c21}
If $p(d+1) \leq 1/e$ then with positive probability no event $A_i$ holds. \npf
\end{coro}

\noindent
{\bf Proof of Theorem \ref{t1}:}\,
Let $a \geq 1$ be any real number, and let $\epsilon > 0$ be small.
Throughout the proof we assume $k$ is sufficiently large as a function of $a$ and $\epsilon$.
Let $k$ be sufficiently large such that there is an integer between
$\frac{k}{(1+\epsilon^2/4)a\ln k}$ and $\frac{k}{(1+\epsilon^2/8)a\ln k}$.
Thus, for some $\epsilon^2/8 \leq \gamma \leq \epsilon^2/4$,
the number $t=\frac{k}{(1+\gamma)a \ln k}$ is an integer.
Now, let $H=(V,E)$ be a hypergraph with $n$ vertices and $\Delta(H) \leq k^a$.
We  will show that there exists an equitable coloring of $H$ with
$\frac{k}{(1+\gamma)a \ln k}- \lceil \sqrt{\gamma}\frac{k}{a \ln k}\rceil > (1-\epsilon)\frac{k}{a\ln k}$
colors.

Assume that we have the set of colors $\{1,\ldots,t\}$.
It will be convenient to deal with the finite set of graphs having $n < 2k\ln k$ separately.
We begin with the general case.

\subsection{The general case: $n > 2k\ln k$}
In the first phase of the proof we color most of the vertices (that is, we obtain a partial coloring) such that
certain very specific properties hold. In the second phase we color the
vertices that were not colored in the first phase and show that we can do it
carefully enough and obtain a proper strong $t$-coloring.
In the third phase we show how to modify our coloring and
obtain an equitable coloring.

\subsubsection{First Phase}
Our goal in this phase is to achieve a partial coloring with the following properties:
\begin{lemma}
\label{l22}
There exists a partial coloring of $H$ with the colors $\{1,\ldots,t\}$ such that the following four conditions hold:
\begin{enumerate}
\item
Every edge contains at least $k\gamma/5$ uncolored vertices.
\item
Every edge has at most $\lceil 10/\gamma \rceil$ colors that do not appear in its vertex set.
\item
Put $z=\lceil k^{1-a\gamma/4} \rceil$. For each $v \in V$, and for each sequence of
$z$ {\bf distinct} colors $c_1,\ldots,c_z$ and for each
sequence of $z$ {\bf distinct} edges containing $v$ denoted $f_1,\ldots,f_z$, at least one $f_i$
has an element colored $c_i$.
\item
Every color appears in at least $n\frac{(1+\gamma/4)a\ln k}{k}$ vertices.
\end{enumerate}
\end{lemma}
{\bf Proof:}\,
We let each vertex $v \in V$ choose a color from $\{1,\ldots,t\}$ randomly.
The probability to choose color $i$ is $p=\frac{(1+\gamma/2)a \ln k}{k}$ for $i=1,\ldots,t$
and the probability of remaining uncolored is, therefore, $q=1-pt=\frac{\gamma}{2(1+\gamma)}$.
For an edge $f$, let $A_f$ denote the event that $f$ contains less than $k\gamma/5$ uncolored vertices.
Let $B_f$ denote the event that $f$ has more than $\lceil 10/\gamma \rceil$ colors missing from its
vertex set. For a vertex $v$, let $C_v$ denote the event that there exist $z$ distinct edges $f_1,\ldots,f_z$
each $f_i$ contains $v$, and there exist $z$ distinct colors $c_1,\ldots,c_z$, such that $c_i$ is
missing from $f_i$ for each $i=1,\ldots,z$. For a color $c$, let $D_c$ denote the event that the
color $c$ appears in less than $n\frac{(1+\gamma/4)a \ln k}{k}$ vertices.
We must show that with positive probability, none of the above $2|E|+|V|+t$ events hold.
The following four claims provide upper bounds for the probabilities of the events $A_f$,$B_f$, $C_v$, $D_c$.
\begin{claim}
\label{c1}
$\Pr[A_f] < \frac{1}{k^{5a}}$.
\end{claim}
{\bf Proof:}\,
Let $X_f$ denote the random variable counting the uncolored elements of $f$.
The expectation of $X_f$ is $E[X_f]=kq$.
Since each vertex chooses its color independently we have by a common Chernoff inequality (cf. \cite{AlSp})
$$
\Pr\left[A_f\right] =\Pr\left[X_f < \frac{k\gamma}{5}\right] \leq \Pr\left[X_f < \frac{kq}{2}\right] = \Pr\left[X_f < \frac{E[X_f]}{2}\right] <
$$
$$
e^{-2(E[X_f]/2)^2/k}=e^{-k^2q^2/(2k)}=e^{-kq^2/2}<< \frac{1}{k^{5a}}. ~~
$$ \npf
\begin{claim}
\label{c2}
$\Pr[B_f] < \frac{1}{k^{5a}}$.
\end{claim}
{\bf Proof:}\,
Fix $s=\lceil 10/\gamma \rceil$ distinct colors. The probability that none of them appear in $f$ is precisely $(1-sp)^k$.
Now,
$$
(1-sp)^k = \left(1-\frac{s(1+\frac{\gamma}{2})a \ln k}{k}\right)^k < \frac{1}{k^{as+as\gamma/2}}.
$$
As there are ${t \choose s} < k^s$ possible sets of $s$ distinct colors
we get that
$$
\Pr[B_f] < {t \choose s}\frac{1}{k^{as+as\gamma/2}} < \frac{1}{k^{as\gamma/2}} \leq \frac{1}{k^{5a}}.
$$ \npf
\begin{claim}
\label{c3}
$\Pr[C_v] < \frac{1}{k^{5a}}$.
\end{claim}
{\bf Proof:}\,
If the degree of $v$ is less than $z$ there is nothing to prove. Otherwise,
fix a set of $z$ distinct colors $\{c_1,\ldots,c_z\}$ and $z$ distinct
edges containing $v$, denoted $\{f_1,\ldots,f_z\}$.
We begin by computing the probability that
for each $i=1,\ldots,z$, $c_i$ does not appear in an element of $f_i$.
Denote this probability by $\rho=\rho(v,f_1,\ldots,f_z,c_1,\ldots,c_z)$ 
Let $\vec{x}=(x_1,\ldots,x_z)$ be a binary vector.
Let $N(\vec{x})$ denote the subset of vertices that belong to the edges $f_i$,
for each coordinate $x_i$ that is positive in $\vec{x}$, and that do not belong to the edge $f_i$
for each coordinate $x_i$ that is zero in $\vec{x}$. This partitions the vertex set into $2^z$ parts.
Thus, if $\vec{x} \neq \vec{x'}$ then $N(\vec{x}) \cap N(\vec{x'})=\emptyset$.
Put $n_{\vec{x}}=|N(\vec{x})|$. Also notice that for each $i=1,\ldots,z$ we have
$(\sum_{\vec{x} \,:\, x_i=1}n_{\vec{x}})=k$.
Let $w_{\vec{x}}$ denote the number of positive coordinates in $\vec{x}$.
Clearly,
$$
\rho=\prod_{\vec{x}} (1-w_{\vec{x}}p)^{n_{\vec{x}}}< \prod_{\vec{x}}e^{-pw_{\vec{x}}n_{\vec{x}}}= e^{-pkz}=\frac{1}{k^{a(1+\gamma/2)z}}.
$$
There are exactly $t^z < (k/\ln k)^z$ ordered sets of $z$ distinct colors.
Thus, the probability that $f_1,\ldots,f_z$ each miss a distinct color
is less than $(k/\ln k)^z/k^{a(1+\gamma/2)z}$. There are at most
${{\lfloor {k^a} \rfloor} \choose z}$ distinct subsets of $z$ edges containing $v$.
This, together with Stirling's formula, gives
$$
\Pr[C_v]  < {{\lfloor {k^a} \rfloor} \choose z}\frac{k^z}{\left(\ln k\right)^zk^{a\left(1+\gamma/2\right)z}} < 
\left(\frac{ek^a}{z}\frac{k}{k^{a\left(1+\gamma/2\right)}\ln k}\right)^z \leq
\left(\frac{e}{k^{a\gamma/4}\ln k}\right)^z << \frac{1}{k^{5a}}.
$$ \npf
\begin{claim}
\label{c4}
$\Pr[D_c] < \frac{1}{e^{n/k}}$.
\end{claim}
{\bf Proof:}\,
Let $X_c$ denote the number of vertices that received the color $c$.
Clearly, $E[X_c] = pn = n \frac{(1+\gamma/2)a \ln k}{k}$.
Put $\beta=n \frac{a \gamma\ln k}{4k}$.
We shall use the Chernoff inequality (cf. \cite{AlSp})
$$
\Pr[X_c - pn < -\beta] < e^{-\beta^2/(2pn)}.
$$
In our case
$$
\Pr[D_c] = \Pr[X_c - pn < -\beta] < e^{-\beta^2/(2pn)}=
e^{-\frac{na\ln k}{k}(\frac{\gamma^2}{32(1+\gamma/2)})} <
$$
$$
e^{-\frac{na\ln k}{k}(\frac{\gamma^2}{33})} = \frac{1}{k^{(n/k)(\gamma^2/33)}} < \frac{1}{e^{n/k}}.
$$ \npf

We now construct a dependency graph for all the events of the form $A_f,B_f,C_v,D_c$
(we refer to the events as ``type A'', ``type B'', ``type C'', and type ``D'' respectively).
Consider an event $A_f$. Let $E(f)$ denote the set of edges of $H$ that are disjoint from $f$.
Let $V(f)$ denote the set of vertices of $H$ that do not appear in any edge that intersects $f$.
Clearly $A_f$ is mutually independent of all the $2|E(f)|+|V(f)|$ events of the form
$A_g$, $B_g$ or $C_v$ which correspond to the elements of $E(f)$ and $V(f)$.
Since there are at most $k^{a+1}$ edges intersecting $f$ and since there are at most
$k^{a+2}$ vertices in these edges, the outdegree in the dependency graph from
$A_f$ to other events of type $A$ is at most $k^{a+1}$. Similarly  the outdegree in the
dependency graph from
$A_f$ to other events of type $B$ is at most $k^{a+1}$, and to events of type $C$ it
is at most $k^{a+2}$. $A_f$ depends on all events of type $D$, so the outdegree is $t$.
This explains the first line of Table \ref{T1} (the dependency table). The other elements in
the table are figured out similarly. Note that events of type $D$ depend on all other
events (the fourth line in Table \ref{T1}).
\begin{table}[h]
\begin{center}
\begin{tabular}{|c|c|c|c|c|}
\hline
source $\backslash$ target & $A_f$ & $B_f$ & $C_v$ & $D_t$ \\ \hline
$A_f$ & $k^{a+1}$ & $k^{a+1}$ & $k^{a+2}$ & $t$\\
$B_f$ & $k^{a+1}$ & $k^{a+1}$ & $k^{a+2}$ & $t$\\
$C_v$ & $k^{2a+1}$ & $k^{2a+1}$& $k^{2a+2}$ & $t$\\
$D_t$ & $|E|$ & $|E|$ & $n$ & $t$\\ \hline
\end{tabular}
\end{center}
\caption{\label{T1}The maximum possible outdegrees in the dependency graph}
\end{table}
In order to apply Lemma \ref{l21} we need to assign a coefficient to each event
in the dependency graph (the coefficients correspond to the $x_i$ in Lemma \ref{l21}).
To each event of type $A$, $B$ or $C$ we assign the coefficient $3/k^{5a}$.
To each event of type $D$ we assign the coefficient $1/e^{n/(2k)}$.
It remains to show that the conditions in Lemma \ref{l21} hold for each event.
Consider events of type $A$. We must show that
$$
\Pr[A_f] < \frac{3}{k^{5a}}\left(1-\frac{3}{k^{5a}}\right)^{k^{a+1}}\left(1-\frac{3}{k^{5a}}\right)^{k^{a+1}}\left (1-\frac{3}{k^{5a}}\right)^{k^{a+2}}\left(1-\frac{1}{e^{n/(2k)}}\right)^{t}.
$$
Indeed, recall that $n >2 k\ln k$ so $(1-1/e^{n/(2k)})^{k-1} > e^{-1}$. Since $t < k-1$ we get, together with Claim \ref{c1},
$$
\frac{3}{k^{5a}}\left(1-\frac{3}{k^{5a}}\right)^{k^{a+1}}\left(1-\frac{3}{k^{5a}}\right)^{k^{a+1}}\left (1-\frac{3}{k^{5a}}\right)^{k^{a+2}}\left(1-\frac{1}{e^{n/(2k)}}\right)^{t} >
$$
$$
\frac{3}{k^{5a}}\left(1-\frac{3}{k^{5a}}\right)^{3k^{a+2}}e^{-1}>
\frac{3}{k^{5a}} \cdot 0.99 \cdot e^{-1}>\frac{1}{k^{5a}} > \Pr[A_f].$$
The analogous inequalities hold for events of type $B$ and $C$ where we use Claim \ref{c2} and Claim \ref{c3} respectively.
Finally, consider events of type $D$. We must show that
$$
\Pr[D_c] < \frac{1}{e^{n/(2k)}}\left(1-\frac{3}{k^{5a}}\right)^{2|E|+n}\left(1-\frac{1}{e^{n/(2k)}}\right)^{t}.
$$
In any $k$-uniform hypergraph, $|E| \leq n\Delta/k$. Thus, in our case, $2|E|+n \leq 3k^{a-1}n$. Using again the fact that$(1-1/e^{n/(2k)})^{k-1} > e^{-1}$
we have, together with Claim \ref{c4},
$$
\frac{1}{e^{n/(2k)}}\left(1-\frac{3}{k^{5a}}\right)^{2|E|+n}\left(1-\frac{1}{e^{n/(2k)}}\right)^{t} > 
$$
$$
\frac{1}{e^{n/(2k)}}\left(1-\frac{3}{k^{5a}}\right)^{3k^{a-1}n}e^{-1}>
\frac{1}{e^{n/(2k)}}\left(1-\frac{3}{k^{5a}}\right)^{(\frac{k^{5a}}{3}-1)\frac{18n}{k^{4a+1}}}e^{-1}>
$$
$$
\frac{1}{e^{n/(2k)}}e^{-\frac{18n}{k^{4a+1}}-1} > \frac{1}{e^{n/(2k)}}\frac{1}{e^{n/(2k)}} = \frac{1}{e^{n/k}}> \Pr[D_c].
$$
According to Lemma \ref{l21}, with positive probability, none of the events in the dependency graph hold. We
have completed the proof of Lemma \ref{l22}.

\subsubsection{Second Phase}
Fix a partial coloring satisfying the four conditions in Lemma \ref{l22}.
For an edge $f$, let $M(f)$ denote the set of missing colors from $f$.
By Lemma \ref{l22} we know that $|M(f)| \leq \lceil 10/\gamma \rceil$.
For a vertex $v$, let $S(v) = \cup_{v \in f} M(f)$.
We claim that $|S(v)| \leq \lceil 10/\gamma \rceil(z-1)\leq 11z/\gamma$.
To see this, notice that if $|S(v)| > \lceil 10/\gamma \rceil(z-1)$ then there must be at least
$z$ distinct edges containing $v$, say, $f_1,\ldots,f_z$ and $z$ distinct colors $c_1,\ldots,c_z$
such that $c_i$ does not appear in $f_i$ for $i=1,\ldots,z$. However, this is impossible
by the third requirement in Lemma \ref{l22}.
In the second phase we only color the vertices that are uncolored after the
first phase. Let $v$ be such a vertex.
We let $v$ choose a random color from $S(v)$ with uniform distribution.
The choices made by distinct vertices are independent (In case $S(v)=\emptyset$ 
we can assign an arbitrary color to $v$).
Let $f \in E$ be any edge, and let $c \in M(f)$. Let $A_{f,c}$ denote
the event that after the second phase, $c$ still does not appear as a color in
a vertex of $f$. Our goal is to show that with positive probability, none of
the events $A_{f,c}$ for $f \in E$ and $c \in M(f)$ hold. This will give a proper strong
$t$-coloring of $H$ (although not necessarily an equitable one).

Let $T(f)$ be the subset of vertices of $f$ that are uncolored after the
first phase. By Lemma \ref{l22} we have $|T(f)| \geq k\gamma/5$. If $c \in M(f)$
we have that for each $u \in T(f)$, the color $c$ appears in $S(u)$.
Hence,
$$
\Pr[A_{f,c}] = \Pi_{u \in T(f)}\left(1-\frac{1}{|S(u)|}\right) \leq \Pi_{u \in T(f)}\left(1-\frac{\gamma}{11z}\right)\leq
$$
$$
\left(1-\frac{\gamma}{11z}\right)^{k\gamma/5}< e^{-\frac{k\gamma^2}{55z}} < e^{-k^{a\gamma/4}\frac{\gamma^2}{110}} << \frac{1}{k^{2a+2}}.
$$
Since each event $A_{f,c}$ is mutually independent of all other events but those that
correspond to edges that intersect $f$, we have that the dependency graph of the events
has maximum outdegree at most $\lceil 10/\gamma \rceil k^{a+1} < k^{a+2}/e-1$.
Since $\frac{1}{k^{a+2}}((k^{a+2}/e-1)+1) = 1/e$ we have, by Corollary \ref{c21},
that with positive probability none of the events of the form $A_{f,c}$ hold.
In particular, there exists a strong $t$-coloring of $H$.

\subsubsection{Third Phase}
Assume the color classes of the strong $t$-coloring obtained after the second phase are $V_1,\ldots,V_t$
where $|V_i| \geq |V_{i+1}|$, $i=1,\ldots,t-1$. By Lemma \ref{l22} we know that
$|V_i| \geq n\frac{(1+\gamma/4)a\ln k}{k}$, $i=1,\ldots,t$. Let $s=\lceil \sqrt{\gamma}k/(a \ln k)\rceil$ and let
$W=V_1 \cup \cdots \cup V_s$. Clearly
$$
n-|W|=|V \setminus W|=|V_{s+1} \cup \cdots \cup V_t|  \geq (t-s) n\frac{(1+\frac{\gamma}{4})a\ln k}{k}=n\left(\frac{1+\frac{\gamma}{4}}{1+\gamma}\right)-\frac{sn(1+\frac{\gamma}{4})a\ln k}{k}.
$$
Hence,
$$
|W| \leq n\left(1-\frac{1+\frac{\gamma}{4}}{1+\gamma}\right)+\frac{sn(1+\frac{\gamma}{4})a\ln k}{k} < \gamma n+\frac{sn(1+\frac{\gamma}{4})a\ln k}{k}.
$$
In particular, $|V_s| \leq |W|/s < \gamma n/s+n(1+\gamma/4)a\ln k/k$. It follows that $\left| |V_i|-|V_j| \right| < \gamma n/s$ for all $s+1 \leq i < j \leq t$.
Hence, it suffices to show that $|W| \geq (t-s)\gamma n/s$ since we can then transfer all the vertices in the color classes $V_1,\ldots,V_s$ to the color
classes $V_{s+1},\ldots,V_t$ such that after the transfer, the $t-s$ remaining classes form an equitable partition (the strong coloring stays proper,
of course). Indeed,
$$
|W| > sn \frac{a\ln k}{k} = s^2n\frac{a\ln k}{sk} \geq n\gamma \frac{k^2}{a^2 (\ln k)^2}\frac{a\ln k}{sk}=n\gamma\frac{k}{sa \ln k}
> n\frac{t\gamma}{s} >(t-s)\frac{n\gamma }{s}.
$$
We have shown how to obtain an equitable coloring with $t-s=\frac{k}{(1+\gamma)a \ln k}- \lceil \sqrt{\gamma}\frac{k}{a \ln k}\rceil > (1-\epsilon)\frac{k}{a\ln k}$
colors.
\subsection{The finite case: $n < 2k\ln k$}
As in the proof for the general case, let each vertex choose a color randomly and independently, each color with
probability $p$ where $p=\frac{(1+\gamma/2)a \ln k}{k}$ for $i=1,\ldots,t$
and the probability of remaining uncolored is $q=1-pt=\frac{\gamma}{2(1+\gamma)}$.
As in the proof of Claim \ref{c1}, the probability that an edge contains less than $k\gamma/5$ uncolored vertices
is less than $1/k^{5a}$. There are $|E| \leq nk^{a}/k \leq 2k^a\ln k$ edges. Hence, the expected number of edges
with less than $k\gamma/5$ edges is less than $1/k^3$. Thus. With probability at least than $1-1/k^3$ all edges
have at least $k\gamma/5$ uncolored vertices. As in the proof of Claim \ref{c4}, the probability that a color appears in less than $na\ln k(1+\gamma/4)/k$ vertices
is less than $\frac{1}{k^{(n/k)(\gamma^2/33)}}$. Unlike Claim \ref{c4} we cannot bound this number from above by
$e^{-n/k}$; instead, since $n \geq k$ (otherwise there are no edges at all), we can bound it with $k^{-\gamma^2/33}$.
Since there are $t < k$ colors, the expected number of colors that appear in less than $na\ln k(1+\gamma/4)/k$ vertices
is less than $k^{1-\gamma^2/33}$. Thus, with probability at least $2/3$ there are less than $3k^{1-\gamma^2/33}$
such colors. Finally, let $X$ count the number of pairs $(e,c)$ where $e \in E$ and $c$ is a color that is missing from
$e$. Clearly, 
$$
E[X]=|E|t(1-p)^k < 2k^a\ln k \cdot k \cdot k^{-a(1+\gamma/2)}=2k^{1-a\gamma/2}\ln k < 2k^{1-\gamma/4} < \frac{k\gamma}{15}.
$$
Hence, with probability at least $2/3$, $X < k\gamma/5$.

We have proved that with probability at least $1-1/k^3-1/3-1/3 > 0$ all the following occur simultaneously:
\begin{enumerate}
\item
All edges have at least $k\gamma/5$ uncolored vertices.
\item
At least $t-3k^{1-\gamma^2/33}$ colors appear each in at least $na\ln k(1+\gamma/4)/k$ vertices.
\item
The number of pairs $(e,c)$ of edges $e$ and colors $c$ such that $c$ is missing from $e$ is
less than $k\gamma/5$.
\end{enumerate}
Fix a partial coloring with all these properties. Trivially we can make it a proper strong coloring by assigning
a color $c$ that is missing from an edge $e$ to one of the uncolored vertices of $e$, and we can do it
greedily to all such $(e,c)$ pairs. We therefore obtain a proper strong $t$-coloring of $H$, where, in addition,
at least $t-3k^{1-\gamma^2/33}$ colors appear each in at least $na\ln k(1+\gamma/4)/k$ vertices.
We can now use the same arguments as in the third phase of the general case and obtain an equitable
coloring. The only difference is that instead of $t$ we only use $t-r$ colors where $r$ is the
number of color classes having less than $na\ln k(1+\gamma/4)/k$ vertices. Thus,
$t-r \geq t-3k^{1-\gamma^2/33}> t(1-\gamma/33)$, and it is easily seen that all
computations in the third phase hold when replacing $t$ with $t(1-\gamma/33)$. \npf

\section{A random hypergraph construction}
Let $a \geq 1$ and let $\epsilon > 0$. Let $n=k^{2a}$. For simplicity we assume $n$ is an integer
in order to ignore floors and ceilings. $k$ will be selected sufficiently large to justify this assumption and
the assumptions that follow. Let $m=(1-\epsilon)k^{3a-1}$ (again, assume $m$ is an integer).
Consider the random $k$-uniform hypergraph on the vertex set $[n]$ with $m$ randomly selected edges
$f_1,\ldots,f_m$. Each edge $f_i$ is chosen uniformly from all ${n \choose k}$
possible edges. The $m$ choices are independent (thus, the same edge can be selected more than once).
The expected degree of a vertex $v$ (including multiplicities) is $mk/n =(1-\epsilon)k^a$. Notice that for $k$ sufficiently
large we have, using a Chernoff inequality, that the degree of $v$ is greater than $k^a$ with probability less than $1/(2k^{2a})=1/(2n)$.
Hence, with probability greater than $0.5$ the maximum degree is at most $k^a$.
Put $t=(1-2\epsilon)na\ln k/k$. Again, we assume $t$ is an integer. We show that with probability greater than $0.5$,
no $t$-subset of vertices is a vertex cover. This proves the existence of hypergraphs $H$
with $\Delta(H) \leq k^a$ and $c(H) \leq (1+o_k(1))k/(a \ln k)$.

Fix $X \subset [n]$ with $|X|=t$. For each edge $f_i$ we have, assuming $k$ is sufficiently large,
$$
\Pr[f_i \cap X = \emptyset] = \frac{(n-t)(n-t-1)\cdots(n-t-k+1)}{n(n-1)\cdots(n-k+1)}> \left(1-\frac{t}{n-k+1}\right)^k >
\left(1-\frac{t}{(1-\epsilon)n}\right)^k=
$$
$$
\left(1-\frac{(1-2\epsilon)a\ln k}{(1-\epsilon)k}\right)^k >  \left(1-\frac{(1-\epsilon)a\ln k}{k}\right)^k > \frac{1}{2}e^{-(1-\epsilon)a\ln k} = \frac{1}{2k^{a(1-\epsilon)}}.
$$
Since each edge is selected independently we have
$$
\Pr[X {\rm~is ~ a ~ vertex ~ cover}] < \left(1-\frac{1}{2k^{a(1-\epsilon)}}\right)^m.
$$
There are ${n \choose t}$ possible choices for $X$. It suffices to show that
$$
{n \choose t}\left(1-\frac{1}{k^{2a(1-\epsilon)}}\right)^m < \frac{1}{2}.
$$
Indeed, for $k$ sufficiently large
$$
{n \choose t}\left(1-\frac{1}{2k^{a(1-\epsilon)}}\right)^m < \left(\frac{en}{t}\right)^t\left(1-\frac{1}{2k^{a(1-\epsilon)}}\right)^{(1-\epsilon)k^{3a-1}}=
$$
$$
\left(\frac{ek}{(1-2\epsilon)a \ln k}\right)^{(1-2\epsilon)k^{2a-1}\ln k}\left(1-\frac{1}{2k^{a(1-\epsilon)}}\right)^{(1-\epsilon)k^{3a-1}}=
$$
$$
\left(\left(\frac{ek}{(1-2\epsilon)a \ln k}\right)^{(1-2\epsilon)\ln k}\left(1-\frac{1}{2k^{a(1-\epsilon)}}\right)^{(1-\epsilon)k^a}\right)^{k^{2a-1}} <
\left(e^{\ln^2 k}e^{-k^{a\epsilon}(1-\epsilon)/2}\right)^{k^{2a-1}} << \frac{1}{2}.
$$ \npf
\section{Concluding remarks}
\begin{itemize}
\item
In the proof of Theorem \ref{t1} we require that $\Delta(H) \leq k^a$ for some fixed $a \geq 1$.
It is possible (although the computations get somewhat more complicated) to prove Theorem \ref{t1}
when $a$ is not necessarily a constant but satisfies $a=a(k) = o(k/\ln k)$. In other words, $\Delta(H)$
is allowed to be any subexponential function of $k$.
\item
The proof of Theorem \ref{t1} is not algorithmic. It is, however, possible to obtain a polynomial time (in the number of vertices
of the hypergraph, and not in its uniformity)
algorithm that yields an equitable partition with $(1-o_k(1))c k/(a \ln k)$ parts where $c$ is a fixed small
constant (depending only on $a$). This can be done by using the method of Beck for the two
coloring of hypergraphs \cite{Be} and generalizing it to more colors. We also need to take care
that the coloring obtained be equitable (Beck's algorithm does not guarantee this). However,
Beck's algorithm can be modified so as to guarantee that all colors use {\em roughly} the same number of colors,
and then we can use the approach from the third phase of our proof to show that by sacrificing only a small fraction of  the colors we can
make the partition equitable using the remaining colors. Notice that the third phase can easily be implemented in polynomial time.
\item
A special case of Theorem \ref{t1} yields an interesting result about graphs. Let $G$ be a $k$-regular graph.
Then, $G$ has an equitable coloring with $(1-o_k(1))(k/\ln k)$ colors such that each color class is a total dominating
set (a total dominating set $D$ is a subset of the vertices that has the property that each vertex $v \in G$ has a neighbor in $D$).
To see this, we can construct a hypergraph $H$ from the graph $G$ as follows. For each vertex $v \in G$
Let $N(v)$ denote the neighborhood of $v$. The vertices of $H$ are those of $G$ and the edges are all the sets $N(v)$.
$H$ is $k$-uniform and $\Delta(H)=k$. Theorem \ref{t1} applied to $H$ gives the desired result about $G$.
\end{itemize}

\end{document}